\documentclass[a4paper,12pt]{article}
\usepackage[cp1251]{inputenc}
\usepackage[ukrainian]{babel}
\usepackage{latexsym,amssymb,amsmath,amsthm,amsfonts,euscript}

\newtheorem{theorem}{Теорема}

\newtheorem{lem}[theorem]{Лема}

\newtheorem{definition}{Означення} 

\title{Аналог теореми Меньшова--Трохимчука для моногенних функцій в тривимірній комутативній алгебрі}

\author{М.В.~Ткачук$^1$ \&\, С.А.~Плакса$^2$\\
Інститут математики НАН України\\ 
$^1$ maxim.v.tkachuk@gmail.com\\[2pt]
$^2$  plaksa62@gmail.com }

\begin{document}

\maketitle

{\bf 1. Вступ.} В алгебрі комплексних чисел $\mathbb {C}$ функція
$F\colon \mathbb {C}\longrightarrow\mathbb {C}$ називається
моногенною в точці $\xi_{0}\in \mathbb {C}$, якщо існує скінченна
границя
\begin{equation} \label{compl-der}
\lim _{\xi\to \xi_{0}}{\frac {F(\xi)-F(\xi_{0})}{\xi-\xi_{0}}}\,.
\end{equation}
При цьому границя (\ref{compl-der})
називається похідною функції $F$ в точці $\xi_{0}$. Функція, яка є
моногенною в усіх точках області $D\subset \mathbb {C}$,
називається голоморфною в цій області (див. \cite{Goursa}).

Встановленню послаблених 
умов голоморфності функцій комплексної змінної присвячені роботи
Х.~Бора \cite{Bohr}, Х.~Радемахера \cite{Rademacher},
Д.Є.~Меньшова \cite{menshov-1,menshov-2,menshov-3}, В.С.~Федорова
\cite{fedorov}, Г.П.~Толстова \cite{Tolstov}, Ю.Ю.~Трохимчука
\cite{Trokhimchuk,zb_trokhinchuk}, Г.Х.~Синдаловського
\cite{Sindalovski}, Є.П.~Долженко \cite{Dolgenko},
 Д.С.~Теляковського \cite{Teliakovski},
М.Т.~Бродович \cite{Brodovich}.

Наведемо тут одну з умов Меньшова, яку, зберігаючи позначення
автора, називають умовою $K'''$, а саме:
кажуть, що {\it функція $F(\xi)$ задовольняє умову $K'''$ в 
точці $\xi_{0}$, 
якщо існує границя (\ref{compl-der}), де $\xi$ належить або
об'єднанню будь-яких двох різних прямих, що перетинаються в точці
$\xi_{0}$, або об'єднанню будь-яких трьох попарно неколінеарних
променів з початком у точці $\xi_{0}$.}

Д.Є.~Меньшов \cite{menshov-1,menshov-2,menshov-3} показав
достатність виконання умови $K'''$ в кожній точці області $D$ (за
винятком не більш ніж зчисленної кількості точок області $D$)
для 
конформності відображення $F$ у випадку, коли $F : D \rightarrow
\mathbb{C}$ --- неперервна однолиста функція. Ю.Ю.~Трохимчук
\cite{Trokhimchuk} зняв умову однолистості функції $F$, довівши
при цьому наступну теорему. \vskip 2mm

\noindent {\bf Теорема Меньшова--Трохимчука.} {\it Якщо функція $F
: D \rightarrow \mathbb{C}$ неперервна в області $D$ і в кожній її
точці, за винятком не більш ніж зчисленної їх кількості,
виконується умова $K'''$, то функція $F$ голоморфна в області $D$. 
} \vskip 2mm

А.В.~Бондар \cite{bondar,bondar-mono} довів аналог цієї теореми
для функцій, заданих в багатовимірному комплексному просторі
$\mathbb{C}^n$, 
при цьому ним доведено, що для голоморфності функції достатньо
неперервності цієї функції та існування і 
рівності похідної Фреше вздовж $2n$ спеціально вибраних напрямів. 
А.В.~Бондар \cite{bondar-mono} і В.І.~Сірик \cite{siryk} довели
також %
для функцій, заданих в 
$\mathbb{C}^n$, аналоги іншої теореми Меньшова--Трохимчука, 
в якій використовується певна умова збереження кутів.
О.С.~Грецький \cite{gretskii} узагальнив згадані тут результати
А.В.~Бондаря на відображення банахових просторів.

Метою даної роботи є послаблення умов моногенності для функцій, що
приймають значення в одній з тривимірних комутативних алгебр над
полем комплексних чисел. При цьому моногенність функції
розуміється як поєднання її неперервності з існуванням похідної
Гато.

\vskip 2mm

{\bf 2. Моногенні функції в тривимірній комутативній гармонічній
алгебрі з двовимірним радикалом.} Розглянемо тривимірну
комутативну асоціативну банахову алгебру
$\mathbb{A}_3$ 
з одиницею $1$ над полем $\mathbb{C}$, базисом якої є трійка $\{1,
\rho, \rho^2\}$, і при цьому виконується рівність $\rho^3 = 0$.
Визначимо евклідову норму елемента алгебри рівністю
$$\|a + b \rho + c \rho^2\|:=\sqrt{|a|^2+|b|^2+|c|^2}\,,\qquad a, b,
c\in \mathbb{C}\,.$$

Алгебра $\mathbb{A}_{3}$ має єдиний максимальний ідеал
$\mathcal{I}:=\{\lambda_{1}\rho+\lambda_{2}\rho^{2}: \lambda_{1},
\lambda_{2}\in\mathbb{C}\}$\, який є також її радикалом.

Оскільки ядром лінійного відображення $f : \mathbb{A}_3
\rightarrow \mathbb{C}$, що визначається рівністю
\begin{equation} \label{mult-func}
f(a + b \rho + c \rho^2) = a\,,
\end{equation}
є максимальний ідеал $\mathcal{I}$, то $f$ є неперервним
мультиплікативним функціоналом (див. \cite[с. 135] {Hil_Filips}).


Зафіксуємо спочатку дійсний тривимірний підпростір  $E_3:=\{\zeta=
xe_1+ye_2+ze_3\,: \,\,x,y,z\in\mathbb{R}\} \subset \mathbb{A}_3$,
де вектори\, $e_1,e_2,e_3$ --- лінійно незалежні над полем дійсних
чисел $\mathbb{R}$, проте, взагалі кажучи, не утворюють базис
алгебри $\mathbb{A}_3$\,. На вибір підпростору $E_3$ накладемо
лише одну вимогу: образом $E_3$ при відображенні $f$ є вся
комплексна площина (див. \cite{Pukh-5,Sh-co}).

Важливими з точки зору застосувань прикладами таких підпросторів є
підпростори, 
побудовані на гармонічних базисах
$\{e_1,e_2,e_3\}$ алгебри $\mathbb{A}_3$\,, що задовольняють
рівність\, $e_1^2+e_2^2+e_3^2=0$ (див.
\cite{melnichenko_plaksa-mono,shpakivskyi_plaksa_umzh}). Існування
гармонічних базисів в комутативній алгебрі є істотною 
передумовою побудови розв'язків тривимірного рівняння Лапласа у
вигляді компонент розкладу диференційовних функцій за векторами
базису (див.
\cite{Ketchum-28,Mel'nichenko75,melnichenko_plaksa-mono}).


Добре відомо, що існують різні типи диференційовності відображень
в лінійних нормованих просторах.  Насамперед, використовуються
сильна дифекренційовність за Фреше і слабка диференційовність за
Гато  (див., наприклад, \cite{Hil_Filips}), при цьому відповідні
похідні Фреше і Гато визначаються як лінійні оператори.  Для
функції, заданої в області скінченновимірної комутативної
асоціативної алгебри, Г.~Шефферс \cite{Scheffers} розглядав
похідну, яка розуміється як функція, визначена в тій самій
області. Узагальнюючи такий підхід 
на випадок довільної комутативної банахової алгебри, Е.Р.~Лорх
\cite{Lorch} ввів сильну похідну функції, яка також розуміється як
функція, визначена в тій же області, що і сама функція.

Функція $\Phi: \Omega \rightarrow \mathbb{A}_3$ називається {\it
диференційовною за Лорхом} в області $\Omega \subset E_3$, якщо
для кожної точки $\zeta \in \Omega$ існує елемент алгебри
$\Phi'_L(\zeta) \in \mathbb{A}_3$ такий, що для кожного
$\varepsilon
> 0$ існує $\delta > 0$ таке, що для всіх $h \in E_3$, для яких
$\|h\| < \delta$, виконується нерівність:
\begin{equation}
\| \Phi(\zeta + h) - \Phi(\zeta) - h \Phi'_L(\zeta) \| \leq \|h\|
\varepsilon.
\end{equation}

Похідна Лорха $\Phi'_L(\zeta)$ є функцією змінної $\zeta$, тобто
$\Phi'_L: \Omega \rightarrow \mathbb{A}_3$.  При цьому
відображення $B_{\zeta}: E_3 \rightarrow \mathbb{A}_3$, задане
рівністю $B_{\zeta} h = h \Phi'_L(\zeta)$, є обмеженим лінійним
оператором. Отже, функція $\Phi$, диференційовна за Лорхом в
області $\Omega$, має похідну Фреше $B_{\zeta}$ в кожній точці
$\zeta \in \Omega$. Обернене твердження загалом не вірне (див.
приклад в монографії \cite[с.~116]{Hil_Filips}).


Використовуючи диференціал Гато,  І.П.~Мельниченко
\cite{Mel'nichenko75} за\-про\-по\-ну\-вав розглядати похідну Гато
також як функцію, визначену в тій же області, що і сама функція.


Якщо для функції $\Phi \colon \Omega\longrightarrow\mathbb{A}_3$,
заданої в області $\Omega\subset E_3$\,, у кожній точці
$\zeta\in\Omega$ існує елемент алгебри
$\Phi_G'(\zeta)\in\mathbb{A}_3$ такий, що
\begin{equation}\label{Gprz}
\lim\limits_{\delta\rightarrow 0+0} \left(\Phi(\zeta+\delta
h)-\Phi(\zeta)\right)\delta^{-1}= h\Phi_G'(\zeta)\quad\forall h\in
E_{3}\,,
\end{equation}
то функцію $\Phi_G' \colon \Omega\longrightarrow\mathbb{A}_3$
будемо називати {\em похідною Гато} функції $\Phi$\,.

Очевидно, що з існування сильної похідної Лорха $\Phi'_L(\zeta)$
випливає існування слабкої похідної Гато $\Phi_G'(\zeta)$ і
рівність $\Phi_L'(\zeta)=\Phi_G'(\zeta)$, проте з існування
похідної Фреше $B_{\zeta}$ не випливає існування похідної
$\Phi_G'(\zeta)$, що демонструє згаданий вище приклад з монографії
\cite[с.~116]{Hil_Filips}.



Розглянемо тепер поняття моногенної функції. 

Функцію $\Phi \colon \Omega \longrightarrow \mathbb{A}_3$
називаємо \textit{моногенною} в області $\Omega\subset E_{3}$\,,
якщо $\Phi$ є неперервною і має похідну Гато в кожній точці
області $\Omega$ (див.
\cite{shpakivskyi_plaksa_umzh,overview_plaksa,Plaksa_UMB}).

Хоча з існування похідної Гато $\Phi_G'(\zeta)$ не випливає
існування похідної Лорха $\Phi'_L(\zeta)$\,, але моногенні функції
$\Phi \colon \Omega \longrightarrow \mathbb{A}_3$ в області
$\Omega\subset E_{3}$ є диференційовними за Лорхом у цій області.
Це випливає з представлення моногенних функцій $\Phi(\zeta)$,
$\zeta\in\Omega$\,, через голоморфні функції комплексної змінної
$f(\zeta)$, встановленого в роботі \cite{shpakivskyi_plaksa_umzh}.

В роботі \cite{Pl-zb17} послабено одну з умов моногенності, а
саме: показано, що за умови існування похідної Гато у функції
$\Phi:\Omega\longrightarrow\mathbb{A}_{3}$ в усіх точках області
$\Omega\subset E_{3}$ неперервність функції $\Phi$ можна замінити
її локальною обмеженістю в області $\Omega$\,.

\vskip 2mm

{\bf 3. Аналог теореми Меньшова--Трохимчука для моногенних функцій
в алгебрі $\mathbb{A}_{3}$.} Введемо деякі позначення.
Перетином радикалу алгебри $\mathbb{A}_3$ з лінійним простором
$E_3$ є множина необоротних елементів, що 
належать $E_3$. Цією множиною є деяка пряма $L:=\{c\, l\, :\, c\in
\mathbb{R}\}$, де через $l\in  E_3$ позначено напрямний вектор
прямої $L$. Прообразом довільної точки $\xi \in \mathbb{C}$ в
$E_3$ при відображенні $f$ є пряма $L^{\zeta}:=\{\zeta+c\, l\, :\,
c\in \mathbb{R}\}$, де $\zeta$
--- деякий елемент із $E_3$ такий, що $\xi=f(\zeta)$. Очевидно, що
пряма $L^{\zeta}$ паралельна прямій $L$.

Зазначимо, що тут i надалі до об’єктів з $E_3$ застосовуються
геометричні поняття (паралельність, опуклість в напрямку прямої
тощо), які,
строго кажучи, 
мають сенс по відношенню до конгруентних прообразів цих об’єктів у
$\mathbb{R}^3$ при взаємно однозначній відповідності
$\zeta=xe_1+ye_2+ze_3$ між елементами $\zeta\in E_3$ і точками
$(x,y,z)\in\mathbb{R}^3$.


Нехай область $\Omega\subset E_3$ є опуклою в напрямку прямої $L$
(область називається {\it опуклою в напрямку прямої}, якщо вона
містить кожен відрізок, який з’єднує дві точки області i
паралельний цій прямій).

Розглянемо наступний гіперкомплексний аналог умови Меньшова $K'''$ 
в алгебрі $\mathbb{A}_3$ для функцій $\Phi : \Omega \rightarrow
\mathbb{A}_3$, визначених в області $\Omega\subset E_3$.



\begin{definition}
Будемо говорити, що функція $\Phi : \Omega \rightarrow
\mathbb{A}_3$ задовольняє умову $K'''_{\mathbb{A}_3, E_3}$ 
в точці $\zeta \in \Omega$, якщо існує елемент $\Phi_*(\zeta)\in
\mathbb{A}_3$ такий, що рівність
\begin{equation}
\lim_{\varepsilon \rightarrow 0,\, \varepsilon \in \mathbb{R}}
\left( \Phi(\zeta + \varepsilon h) - \Phi(\zeta) \right)
\varepsilon^{-1} = h \Phi_*(\zeta) 
\label{eq:deriv}
\end{equation}
виконується для трьох векторів $h$, а саме: 
векторів $h_1, h_2, l$, що утворюють базис в просторі $E_3$.
\end{definition}

Зауважимо, що у випадку, коли функція $\Phi : \Omega \rightarrow
\mathbb{A}_3$ задовольняє умову $K'''_{\mathbb{A}_3, E_3}$ в
різних точках області $\Omega \subset E_3$, набір векторів $h_1,
h_2$ може бути різним в різних точках цієї області.


\begin{lem}
\label{monom_K} Нехай область $\Omega \subset E_3$ є опуклою в
напрямку прямої $L$ і неперервна в $\Omega$ функція $\Phi : \Omega
\rightarrow \mathbb{A}_3$ має вигляд $\Phi(\zeta) = \rho^2
\Phi_2(\zeta)$, де $\Phi_2(\zeta) \in \mathbb{C}$, і задовольняє
умову $K'''_{\mathbb{A}_3, E_3}$ в усіх точках $\zeta \in \Omega$,
крім не більш ніж зчисленної множини точок. Тоді $\Phi_2(\zeta) =
F_2(f(\zeta))$, де $F_2 : D \rightarrow \mathbb{C}$ --- голоморфна
функція в області $D$, 
яка є образом області $\Omega$ при відображенні $f$.
\end{lem}

\begin{proof} Нехай $\zeta \in \Omega$ --- довільна точка, в якій функція $\Phi$ задовольняє умову $K'''_{\mathbb{A}_3,
E_3}$. Запишемо рівність (\ref{eq:deriv}) для функції $\Phi(\zeta)
= \rho^2 \Phi_2(\zeta)$:
\begin{equation}
\lim_{\varepsilon \rightarrow 0,\, \varepsilon \in \mathbb{R}}
\rho^2 \left( \Phi_2(\zeta + \varepsilon h) - \Phi_2(\zeta)
\right) \varepsilon^{-1} = h \Phi_*(\zeta) \label{eq:deriv_monom}
\end{equation}
і зазначимо, що вона виконується при $h\in\{h_1, h_2, l\}$.

Підставимо $h = h_1$ у рівність (\ref{eq:deriv_monom}) і з
урахуванням того, що $h_1$ є оборотним елементом алгебри
$\mathbb{A}_3$, отримаємо
\begin{equation} \label{eq:deriv_monom_1}
\Phi_*(\zeta) = \rho^2\, h_1^{-1} \lim_{\varepsilon \rightarrow
0,\, \varepsilon \in \mathbb{R}} \left( \Phi_2(\zeta + \varepsilon
h_1) - \Phi_2(\zeta) \right) \varepsilon^{-1} =: \rho^2\,
\Psi(\zeta).
\end{equation}

Після підстановки виразу (\ref{eq:deriv_monom_1}) для $\Phi_*$ в
рівність (\ref{eq:deriv_monom}) вона набуде вигляду
\begin{equation}
\lim_{\varepsilon \rightarrow 0,\, \varepsilon \in \mathbb{R}}
\rho^2 \left( \Phi_2(\zeta + \varepsilon h) - \Phi_2(\zeta)
\right) \varepsilon^{-1} = h \rho^2 \Psi(\zeta).
\label{eq:deriv_monom_2}
\end{equation}

Тепер після підстановки в (\ref{eq:deriv_monom_2}) значення $h =
l$ отримаємо нуль в правій частині рівності
(\ref{eq:deriv_monom_2}). Звідси випливає, що похідна функції
$\Phi_2$ вздовж прямої $L^{\zeta}$ дорівнює нулю всюди, крім не
більш ніж зчисленної множини точок. Тоді за теоремою 9 з
монографії Ю.Ю.~Трохимчука \cite[с.
103]{zb_trokhinchuk} функція $\Phi_2$ є сталою на 
перетині області $\Omega$ з прямою $L^{\zeta}$ (при цьому 
перетини області $\Omega$ з усіма прямими $L^{\zeta}$ є зв'язними
внаслідок опуклості області $\Omega$ в напрямку прямої $L$). Отже,
функція $\Phi_2$ може бути представлена у вигляді $\Phi_2(\zeta) =
F_2(f(\zeta))$, де $F_2: D \rightarrow \mathbb{C}$ --- деяка
неперервна в області $D$ функція.

Доведемо, що функція $F_2$ голоморфна в області $D$.

Спочатку зазначимо, що наслідком означення (\ref{mult-func})
функціонала $f$ є рівність
$$ \rho^2\, h \Psi(\zeta)=\rho^2\, f(h) f(\Psi(\zeta)).$$
Тому, позначаючи при цьому $\xi:=f(\zeta)$, переписуємо рівність
(\ref{eq:deriv_monom_2}) у вигляді
\begin{equation}
\rho^2\,\lim_{\varepsilon \rightarrow 0,\, \varepsilon \in
\mathbb{R}}
 \left( F_2(\xi + \varepsilon f(h)) - F_2(\xi) \right)
\varepsilon^{-1} = \rho^2\, f(h) f(\Psi(\zeta)).
\label{eq:deriv_monom_3}
\end{equation}

Оскільки вирази біля $\rho^2$ в обох частинах рівності
(\ref{eq:deriv_monom_3}) приймають комплексні значення, то з
єдиності розкладу елемента алгебри за базисом випливає рівність
$$
\lim_{\varepsilon \rightarrow 0,\, \varepsilon \in \mathbb{R}}
\left( F_2(\xi + \varepsilon f(h)) - F_2(\xi) \right)
\varepsilon^{-1} = f(h) f(\Psi(\zeta)),
$$
яка виконується при $h \in \{h_1, h_2\}$.

Звідси випливають рівності
\begin{multline*}
f(\Psi(\zeta)) = \lim_{\varepsilon \rightarrow 0,\, \varepsilon
\in \mathbb{R}} \left( F_2(\xi + \varepsilon t_1) - F_2(\xi)
\right)
(\varepsilon t_1)^{-1} =\\
=\lim_{\varepsilon \rightarrow 0,\, \varepsilon \in \mathbb{R}}
\left( F_2(\xi + \varepsilon t_2) - F_2(\xi) \right) (\varepsilon
t_2)^{-1},
\end{multline*}
де\,\, $t_1 := f(h_1)$, $t_2 := f(h_2)$.


Отже, в кожній точці області $D$, за винятком не більш ніж
зчисленної їх кількості, існують похідні функції $F_2$ вздовж двох
різних прямих і ці похідні рівні, а це означає, що неперервна
функція $F_2$ задовольняє умову $K'''$ Меньшова.  Тоді з теореми
Меньшова--Трохимчука
випливає голоморфність функції $F_2$ в області $D$.
\end{proof}



Кожен елемент $a + b \rho + c \rho^2$, $a, b, c\in \mathbb{C}$, за
умови $a\ne 0$ має обернений елемент, розклад якого за 
базисом $\{1,\rho,\rho^2\}$ визначається рівністю
$$(a + b \rho + c \rho^2)^{-1}=\frac{1}{a}-\frac{b}{a^2}\,\rho +\left(\frac{b^2}{a^3}-\frac{c}{a^2}\right)\rho^2\,.$$
Використовуючи цей розклад, легко виписати розклад за 
базисом $\{1,\rho,\rho^2\}$ головного продовження голоморфної
функції $F : D\rightarrow\mathbb{C}$ в область $\Pi:=\{\zeta\in
E_3 : f(\zeta)\in D\}$, 
яка очевидно є нескінченним циліндром, твірні якого паралельні
прямій $L$:
\begin{multline}\label{gol-prod}
\frac{1}{2\pi
i}\int\limits_{\gamma}F(t)(t-\zeta)^{-1}\,dt= F(f(\zeta)) 
+(a_1x+b_1y+c_1z)F'(f(\zeta))\,\rho+\\ 
+\biggl((a_2x+b_2y+c_2z)F'(f(\zeta))+\frac{(a_1x+b_1y+c_1z)^2}{2}F''(f(\zeta))\biggr)\,\rho^{2}\\ 
\forall\zeta=xe_{1}+ye_{2}+ze_{3}\in\Pi\,,
\end{multline}
 де $i$ --- уявна комплексна одиниця, замкнена жорданова спрямлювана крива $\gamma$
лежить в області $D$ і охоплює точку $f(\zeta)=a_0x+b_0y+c_0z$, а
комплексні сталі $a_k, b_k, c_k$ при $k=0,1,2$, --- це коефіцієнти
з розкладів елементів $e_{1}, e_{2}, e_{3}$ за базисом
$\{1,\rho,\rho^2\}$:
$$\begin{array}{l}
    e_1=a_0+a_1\rho+a_2\rho^2,\\
    e_2=b_0+b_1\rho+b_2\rho^2,\\
    e_1=c_0+c_1\rho+c_2\rho^2.
   \end{array} $$
Розклад (\ref{gol-prod}) узагальнює аналогічний розклад, отриманий
в теоремі 1.7 з \cite{melnichenko_plaksa-mono} при додатковому
припущенні, що $e_1=1$.


\begin{lem}
\label{monom_K-1} Нехай область $\Omega \subset E_3$ є опуклою в
напрямку прямої $L$, функція $\Phi : \Omega \rightarrow
\mathbb{A}_3$ є неперервною в $\Omega$ і задовольняє умову
$K'''_{\mathbb{A}_3, E_3}$ в усіх точках $\zeta \in \Omega$\,,
крім не більш ніж зчисленної множини точок. Тоді при всіх $\zeta
\in \Omega$ справедливе представлення
\begin{equation}\label{predst}
\Phi(\zeta) = \frac{1}{2\pi i}\int\limits_{\gamma}
\Big(F_0(\xi)+F_{1}(\xi)\rho+F_{2}(\xi)\,\rho^2\Big) (\xi -
\zeta)^{-1}\, d \xi\,,
\end{equation}
де $F_0, F_1, F_2$ --- деякі функції, голоморфні в області $D$,
яка є образом області $\Omega$ при відображенні $f$.
\end{lem}

\begin{proof} При $\zeta \in \Omega$ розглянемо розклад $\Phi(\zeta)$ за базисом $\{1,\rho,\rho^2\}$:
$$\Phi(\zeta) = \Phi_0(\zeta) + \Phi_1(\zeta) \rho + \Phi_2(\zeta) \rho^2.$$
 Функція $\rho^2 \Phi(\zeta) = \rho^2
\Phi_0(\zeta)$ є неперервною в $\Omega$ і задовольняє умову
$K'''_{\mathbb{A}_3, E_3}$ в усіх точках $\zeta \in \Omega$, крім
не більш ніж зчисленної множини точок. Тоді з леми \ref{monom_K}
випливає, що $\Phi_0(\zeta) = F_0(f(\zeta))$, де $F_0$ ---
голоморфна функція в області $D$, яка є образом області $\Omega$
при відображенні $f$.

Як випливає з рівності (\ref{gol-prod}), перші компоненти в
розкладах за базисом $\{1,\rho,\rho^2\}$ функцій $\Phi(\zeta)$ і
$\frac{1}{2\pi i}\int_{\gamma} F_0(\xi) (\xi - \zeta)^{-1}\, d
\xi$ співпадають в області $\Omega$. Тому справедлива рівність
\begin{equation}\label{rivn1}
\Phi(\zeta) - \frac{1}{2\pi i}\int\limits_{\gamma} F_0(\xi) (\xi -
\zeta)^{-1}\, d \xi = \Phi_{11}(\zeta)\, \rho + \Phi_{12}(\zeta)\,
\rho^2 \quad \forall\zeta \in \Omega,
\end{equation}
 де
$\Phi_{11}$, $\Phi_{12}$ --- деякі комплекснозначні неперервні в
$\Omega$ функції.

Тоді функція $\rho (\Phi_{11}(\zeta) \rho + \Phi_{12}(\zeta)
\rho^2) = \rho^2 \Phi_{11}(\zeta)$  є неперервною в $\Omega$ і
задовольняє умову $K'''_{\mathbb{A}_3, E_3}$ в усіх точках $\zeta
\in \Omega$, крім не більш ніж зчисленної множини точок. Отже, за
лемою \ref{monom_K} маємо $\Phi_{11}(\zeta) = F_1(f(\zeta))$, де
$F_1$
--- голоморфна функція в області $D$.

Далі так, як і при доведенні рівності (\ref{rivn1}), отримуємо
рівність
\begin{equation}\label{rivn2}
\Phi_{11}(\zeta)\, \rho + \Phi_{12}(\zeta)\, \rho^2 -
\rho\,\frac{1}{2\pi i}\int\limits_{\gamma} F_1(\xi) (\xi -
\zeta)^{-1}\, d \xi = \Phi_{22}(\zeta)\, \rho^2 \quad \forall\zeta
\in \Omega,
\end{equation}
 де $\Phi_{22}$ --- деяка комплекснозначна неперервна в
$\Omega$ функція.

Як наслідок рівностей (\ref{rivn1}), (\ref{rivn2}), маємо рівність
\begin{multline}\label{rivn3}
\Phi(\zeta) - \frac{1}{2\pi i}\int\limits_{\gamma} F_0(\xi) (\xi -
\zeta)^{-1}\, d \xi -\\
-\rho\, \frac{1}{2\pi
i}\int\limits_{\gamma} F_{1}(\xi) (\xi - \zeta)^{-1}\, d \xi =
\Phi_{22}(\zeta) \rho^2 \quad \forall\zeta \in \Omega.
\end{multline}

Тепер, спираючись на лему \ref{monom_K}, приходимо до рівності
$\Phi_{22}(\zeta) = F_2(f(\zeta))$, де $F_2$
--- голоморфна функція в області $D$. Тому справедливими є також
рівності
\begin{equation}\label{rivn4}
\rho^2\,\Phi_{22}(\zeta)=\rho^2\,F_2(f(\zeta))  = \rho^2\,
\frac{1}{2\pi i}\int\limits_{\gamma} F_{2}(\xi) (\xi -
\zeta)^{-1}\, d\xi \quad \forall\zeta \in \Omega\,.
\end{equation}

Нарешті, як наслідок рівностей (\ref{rivn3}), (\ref{rivn4}),
отримуємо представлення (\ref{predst}).
\end{proof}

Основним результатом пункту 3 є наступне твердження.

\begin{theorem} \label{theor-1}
Нехай область $\Omega \subset E_3$ є опуклою в напрямку прямої
$L$, функція $\Phi : \Omega \rightarrow \mathbb{A}_3$ є
неперервною в $\Omega$ і задовольняє умову $K'''_{\mathbb{A}_3,
E_3}$ в усіх точках $\zeta \in \Omega$\,, крім не більш ніж
зчисленної множини точок. Тоді:

1) функція $\Phi$ є моногенною в області $\Omega$;

2) функція $\Phi$ продовжується до функції, моногенної в області
$\Pi$. Таке продовження єдине і задається рівністю (\ref{predst})
при всіх $\zeta\in\Pi$;

3) моногенне продовження (\ref{predst}) функції $\Phi$  є
диференційовним за Лорхом в області $\Pi$.
\end{theorem}

Усі твердження теореми \ref{theor-1} є очевидними наслідками
пред\-став\-лен\-ня (\ref{predst}).


\bigskip



\begin{thebibliography}{9}

\bibitem{Goursa}
Goursat~E. Cours d'analyse mathematique.~---  Paris:
Gauthier--Villars, 1910.~--- Vol.~2.

\bibitem{Bohr}
Bohr~H. \"Uber streckentreue und konforme Abbildung // Math.
Zeitschr., 1918, {\bf 1}, 403--420.

\bibitem{Rademacher}
Rademacher~H. \"Uber streckentreue und winkeltreue Abbildung //
Math. Zeitschr., 1919, {\bf 4}, 131--138.

\bibitem{menshov-1}
Menchov~D. Sur les differentielles totales des fonctions
univalentes // Math. Ann. {\bf 105} (1931), 75–85.

\bibitem{menshov-2}
Menchov~D. Sur les fonctions monogenes // Bull. Soc. math. France.
1931. {\bf 59}. С 141-182.


\bibitem{menshov-3}
Menchov~D. Les conditions de monogeneite // Act. Sci. et Ind.
1936. № 329. Paris

\bibitem{fedorov}
Федоров~В.С. О моногенных функциях // Мат. сб., 1935, {\bf 42}, №
4, 485–500

\bibitem{Tolstov}
Толстов~Г.П. О криволинейном и повторном интеграле~// Труды Мат.
ин-та АН СССР.
--- 1950. --- {\bf 35}. ---  C. 3 --- 101.

\bibitem{Trokhimchuk}
Трохимчук~Ю.Ю. Непрерывные отображения и условия моногенности.
--- Москва: Физматиз, 1963. --- 212 с.

\bibitem{zb_trokhinchuk}
Трохимчук~Ю.Ю. Дифференциирование, внутренние отображения и
критерии аналитичности. - Киев : Ін-т математики НАН Украины,
2007. – 539 с. – (Математика та її застосування) (Праці / Ін-т
математики НАН України ; т. 70).

\bibitem{Sindalovski}
Синдаловский~Г.Х. О дифференцируемости и аналитичности однолистных
отображений // Докл. АН СССР. 1979. {\bf 249}. № 6. С. 1325-27

\bibitem{Dolgenko}
Долженко~Е.П. Работы Д.Е.~Меньшова по теории аналитических функций
и современное состояние теории моногенности // УМН, 1992, {\bf
47}, № 5, 67--96.

\bibitem{Teliakovski}
Теляковский~Д.С. Об ослаблении условия асимптотической
моногенности // Мат. заметки, 1996, {\bf 60}, № 6, 902--911.

\bibitem{Brodovich}
Бродович~М.Т. Об отображениях пространственной области,
сохраняющих углы и растяжения вдоль системы лучей // Сиб. мат.
журн. - 1997. - {\bf 38}, № 2. - С. 260-262.

\bibitem{bondar}
Бондарь А.В. Многомерное обобщение одной теоремы Д. Е. Меньшова //
Укр. мат. журн. - 1978. - {\bf 30}, № 4. - С. 435–443.

\bibitem{bondar-mono}
Бондарь А.В. Локальные геометрические характеристики голоморфных
отображений. -- Киев: Наук. думка, 1992. -- 220 с.

\bibitem{siryk}
Сірик В.І. Некоторые критерии голоморфности непрерывных
отображений. - Укр. мат. журн. - 1985. - 37, № 6. - С. 751–756

\bibitem{gretskii}
Грецький О.С. Про C-диференційовиість відображень банахових
просторів. - Укр. мат. журн. - 1994. - 46, № 10. - С. 1336–1342

\bibitem{Hil_Filips}
\emph{Hille~E., Phillips~R.\,S.} Functional Analysis and
Semi-Groups.~--- Providence, R.I.:  Amer. Math. Soc., 1957.



\bibitem{Pukh-5}
Plaksa~S.A., Pukhtaievych~R.P. \emph{Monogenic functions in a
finite-dimensional semi-simple commutative algebra}~// An. St.
Univ. Ovidius Constanta, \textbf{22}, (2014),  N 1, 221--235.

\bibitem{Sh-co}
Shpakivskyi~V. \emph{Constructive description of monogenic
functions in a finite-dimensional commutative associative
algebra}~// Adv. Pure Appl. Math., \textbf{7}, (2016), N 1,
63–-75.

\bibitem{melnichenko_plaksa-mono}
Мельниченко~И.П., Плакса~С.А. Коммутативные алгебры и
пространственные поля.--- Киев: Ин-т математики НАН Украины, 2008.
– 230 с. 

\bibitem{shpakivskyi_plaksa_umzh}
Plaksa~S.A., Shpakovskii~V.S. Constructive description of
monogenic functions in a harmonic algebra of the third rank //
Ukr. Math. J., 62 (2011), no. 8, 1251–1266.

\bibitem{Ketchum-28}
Ketchum~P.W. \emph{Analytic functions of hypercomplex
variables}~// Trans. Amer. Math. Soc., \textbf{30}, (1928),
641--667.

\bibitem{Mel'nichenko75}
Mel'nichenko~I.P. \emph{The representation of harmonic mappings by
monogenic functions}~// Ukr. Math. J., \textbf{27}, (1975), N 5,
499--505.

\bibitem{Scheffers}
Scheffers~G.  \emph{Verallgemeinerung der grundlagen der
gew\"ohnlich complexen fuktionen, I, II}~// Ber. Verh. Sachs.
Akad. Wiss. Leipzig Mat.-Phys. Kl., \textbf{45}, (1893), 828-–848;
\textbf{46}, (1894), 120-–134.

\bibitem{Lorch}
Lorch~E.R.  The theory of analytic function in normed abelian
vector rings // Trans. Amer. Math. Soc., 54 (1943), 414–425.

\bibitem{overview_plaksa}
Plaksa~S.A. Commutative algebras associated with classic equations
of mathematical physics // Advances in Applied Analysis, Trends in
Mathematics, Springer, Basel (2012), 177–223.

\bibitem{Plaksa_UMB}
Plaksa~S.A. Monogenic functions in commutative algebras associated
with classical equations of mathematical physics // Укр. мат.
вісник, 2018, {\bf 15}, № 4, 543–575.

\bibitem{Pl-zb17}
Plaksa~S.A. \emph{On differentiable and monogenic functions in a
harmonic algebra}~// Zb. Pr. Inst. Mat. NAN Ukr., \textbf{14},
(2017),  N 1, 210--221.


\end{thebibliography}
\end{document}